\documentclass[12pt]{article}

\usepackage{amssymb}
\usepackage{amsmath}
\usepackage{graphicx}
\usepackage{hyperref}

\numberwithin{equation}{section}

\begin{document}

\newcommand{\bd}{\begin{displaymath}}
\newcommand{\ed}{\end{displaymath}}
\newcommand{\ds}{\displaystyle}
\newcommand{\bp}{\underline{\bf Proof}:\ }
\newcommand{\ep}{{\hfill $\Box$}\\ }
\newcommand{\be}{\begin{equation}}
\newcommand{\ee}{\end{equation}}
\newcommand{\ba}{\begin{array}}
\newcommand{\ea}{\end{array}}
\newcommand{\bea}{\begin{eqnarray}}
\newcommand{\eea}{\end{eqnarray}}
\newcommand{\nt}{\noindent}

\newtheorem{0}{DEFINITION}[section]
\newtheorem{1}{LEMMA}[section]
\newtheorem{2}{THEOREM}[section]
\newtheorem{3}{COROLLARY}[section]
\newtheorem{4}{PROPOSITION}[section]
\newtheorem{5}{REMARK}[section]
\newtheorem{6}{EXAMPLE}[section]
\newtheorem{7}{ALGORITHM}[section]
\newtheorem{8}{CONJECTURE}[section]

\title{A Graph Theory Assisted Method for Determining Limiting Probability Distributions of Higher Order Markov Chains}
\author{
Jianhong Xu\thanks{School of Mathematical and Statistical Sciences, Southern Illinois University Carbondale, Carbondale, IL 62901, USA. Email: \texttt{jhxu@siu.edu}} 
}

\maketitle

\begin{abstract}  
This paper presents a graph theory assisted method that checks a sufficient condition for the existence of the limiting probability distribution of a higher order Markov chain and identifies this distribution when the condition is satisfied. This method is more effective as compared to the direct verification of the condition since the latter hinges wholly on computations of the full spectrum on a large-size transition matrix. Even though the graph theoretic tools employed here are well known, their utilization in the context of a higher order Markov chain is new. This method is also implemented as an easy-to-use MATLAB function {\tt gtam}. Several illustrative examples are provided too. 
\end{abstract}

\nt {\bf Keywords}: higher order Markov chains, limiting probability distributions, transition tensors, digraphs, strongly connected components, aperiodicity

\nt {\bf AMS Subject Classification}: 05C50, 05C85, 15A18, 15B51, 60J10, 60J22 

\section{Introduction}
\label{intro}
\setcounter{equation}{0}

An $(m-1)$th order, where the integer $m \ge 2$, Markov chain on state space $S=\{1, 2, \ldots, n\}$ refers to a stochastic process $X=\{X_t : t=1, 2, \ldots\}$ such that 
\begin{gather*}
\Pr(X_{t+1}=i_1 | X_t=i_2, \ldots, X_{t-m+2}=i_m, \ldots, X_1=i_{t+1})\\
=\Pr(X_{t+1}=i_1 | X_t=i_2, \ldots, X_{t-m+2}=i_m)
\end{gather*}
for any $t \ge m-1$ and $i_1, \ldots, i_m, \ldots, i_{t+1} \in S$. Such a chain is called homogeneous if the above probability is independent of the parameter $t$,  which is  regarded as time. Accordingly, this probability can be simply denoted as 
 $$\Pr(X_{t+1}=i_1 | X_t=i_2, \ldots, X_{t-m+2}=i_m)=p_{i_1i_2\ldots i_m}$$
 and called the transition probability from states $(i_2, \ldots, i_m)$ to $i_1$. Throughout this paper, we shall always assume such homogeneity. The chain is said to be higher order when $m \ge 3$, which will be the focus of this paper. Nevertheless, our main results remain valid in the first order case.
 
The $m$th order, $n$ dimensional tensor ${\cal P}=[p_{i_1i_2\ldots i_m}]$ consisting of all the transition probabilities of the chain is called its transition tensor.
When the chain is of first order, $\cal P$ is more commonly known as the transition matrix, and is also denoted as $P$ when distinguishing it from the higher order case is needed.

The probability distribution of the chain at time $t$ is defined as  
\be
\label{pd}
x_t=[\Pr(X_t=1) ~\Pr(X_t=2) ~\ldots ~\Pr(X_t=n)]^T.
\ee
The following notion has been introduced in \cite{HX1}.

\begin{0}
For an $(m-1)$th order, $n$-state Markov chain, let $x_t$ be given as in (\ref{pd}). If $\lim_{t \to \infty}x_t$ exists and if this limit does not depend on the initial probability distributions $x_1, x_2, \ldots, x_{m-1}$, then $\pi=\lim_{t \to \infty}x_t$ is called the limiting probability distribution of the chain.
\end{0}

Such a limiting probability distribution manifests the probabilities of the chain to be eventually in each of its states. Hence it serves as an important characterization for the long-term behavior of the chain and naturally raises a key question regarding conditions for its existence.

For a first order Markov chain, a well-known condition that ensures the existence of its limiting probability distribution is the regularity of its transition matrix $P$ \cite{KS, Mey}, i.e., there exists integer $k \ge 1$ such that $P^k > 0$. The chain is called regular as well in this case. Under such regularity condition, $\lim_{k \to \infty}P^k =\pi e^T$, where $\pi > 0$ is the limiting probability distribution and $e$ is the vector of all ones. In addition, $\pi$ is the unique (right) eigenvector of $P$, normalized in the sense that $\pi > 0$ and $\|\pi\|_1=1$, corresponding to the dominant eigenvalue $\lambda=1$. In matrix theory, this type of regularity is also called the primitivity of $P$.

An earlier attempt to extend the above result to a second order Markov chain can be found in \cite{Vla}. A more recent attempt in the same direction but involves a general higher order Markov chain has been made in \cite{HX1}. In both works, in particular, the regularity of the transition tensor has been claimed to lead to the existence of the limiting probability distribution of the chain. Such a result, however, turns out to be incorrect \cite{HX2}.

To address the above problem, a new sufficient condition has been proposed in \cite{HX2} for a general higher order Markov chain. It is formulated using the eigenvalues and eigenvectors of the transition matrix of a relevant first order chain on a much larger state space with $N=n^{m-1}$ states. Hence directly verifying this condition requires computing the full spectrum of that matrix, which typically has time complexity $O(N^3)$ \cite{GV}. Such computations can become prohibitively expensive even for moderate values of $m$ and $n$.

In view of this situation, our primary goals in this paper are to describe an effective graph theory assisted approach to checking the aforementioned sufficient condition in \cite{HX2} for a higher order Markov chain and to present its software implementation in MATLAB. A few illustrative examples will also be provided. Although the graph theoretic tools employed in this paper are well established, their application to a higher order Markov chain appears to be new as far as we can see.

\section{More Background Material}
\label{back}
\setcounter{equation}{0}

Recall that $\cal P$ is the transition tensor of an $(m-1)$th order, $n$ dimensional Markov chain $X=\{X_t : t=1, 2, \ldots\}$ on the state space $S$. This $\cal P$ is also called the one-step transition tensor. More generally, for integer $k \ge 1$, the $k$-step transition probabilities of the chain, i.e., 
$$\Pr(X_{t+k}=i_1 | X_t=i_2, \ldots, X_{t-m+2}=i_m)=p^{(k)}_{i_1i_2\ldots i_m},$$ form the $k$-step transition tensor ${\cal P}^{(k)}=[p^{(k)}_{i_1i_2\ldots i_m}]$. In particular, ${\cal P}^{(1)}=\cal P$. Similar to the first order case, $\cal P$ (and the chain $X$ as well) is called regular if ${\cal P}^{(k)} > 0$ for some $k \ge 1$. 

By convention, ${\cal P}^{(0)}=[p^{(0)}_{i_1i_2\ldots i_m}]$ equals to the $m$th order, $n$ dimensional identity tensor ${\cal I}=[\delta_{i_1i_2\ldots i_m}]$, where for any $i_1, i_2, \ldots, i_m \in S$,
$$\delta_{i_1i_2\ldots i_m}=\left\{\begin{array}{cl}
1, & i_1=i_2;\\
0, & {\rm otherwise.}
\end{array}\right.$$
To determine ${\cal P}^{(k)}$, it is known \cite{HWX} that 
$$p^{(k+1)}_{i_1i_2\ldots i_m}=\sum_{j \in S}p^{(k)}_{i_1ji_2\ldots i_{m-1}}p_{ji_2\ldots i_m}, ~k=0, 1, 2, \ldots,$$
which can also be expressed using the tensor $\boxtimes$ product \cite{HX1, Xu} as 
$${\cal P}^{(k+1)}={\cal P}^{(k)} \boxtimes {\cal P}, ~k=0, 1, 2, \ldots.$$

From now on, $N=n^{m-1}$. Denote by $P^{(k)} \in \mathbb R^{n \times N}$ the mode-$1$ matricization of ${\cal P}^{(k)}$, consisting of all the frontal slices ${\cal P}^{(k)}(:,:,i_3, \ldots, i_m)$, arranged side by side in the linear indexing order of $(i_3, \ldots, i_m)$ \cite{KB}. Especially, $P^{(1)}$ is also denoted simply as $P$.

It is well known \cite{Doo} that the vector-valued stochastic process 
$$Y=\{Y_t : t=m-1, m, \ldots\},$$
where $Y_t=[X_t ~X_{t-1} ~\ldots ~X_{t-m+2}]^T$, is a first order Markov chain with state space $T=\{i_1i_2\ldots i_{m-1} : i_1, i_2, \ldots, i_{m-1} \in S\}$. For brevity, the states in $T$ are denoted here in multi-index form, arranged by the linear indexing order of $(i_1, i_2, \ldots, i_{m-1})$. The notation $Y_t=i_1i_2\ldots i_{m-1}$ means $X_t=i_1, X_{t-1}=i_2, \ldots, X_{t-m+2}=i_{m-1}$.
This chain $Y$ is called the reduced first order chain associated with $X$. Its transition matrix can be denoted in multi-index form as $Q=[q_{i_1i_2\ldots i_{m-1}, j_2j_3\ldots j_m}] \in \mathbb R^{N \times N}$. Then, the entries of $Q$ can be expressed conveniently as \cite{HX1}
$$q_{i_1i_2\ldots i_{m-1}, j_2j_3\ldots j_m}=\left\{\begin{array}{cl}
p_{i_1i_2\ldots i_{m-1}j_m}, & i_\ell=j_\ell, ~\ell=2, 3, \ldots, m-1;\\
0, & {\rm otherwise.}
\end{array}\right.$$
Clearly, $Q$ is column stochastic.

The probability distribution $y_t \in \mathbb R^N$ of the chain $Y$ at time $t$ is analogous to (\ref{pd}). Specifically, the $i_1i_2\ldots i_{m-1}$th entry of $y_t$ is given by $\Pr(Y_t=i_1i_2\ldots i_{m-1})$. It is known \cite{LN} that for $t \ge m-1$, 
\be
\label{xt}
x_{t+1}=Py_t
\ee
and
\be
\label{yt}
y_{t+1}=Qy_t.
\ee

The following result establishes a sufficient condition for the existence of the limiting probability distribution of a higher order Markov chain \cite{HX2}.

\begin{2}
\label{suff}
Let $\cal P$ be the transition tensor of an $(m-1)$th order, $n$-state Markov chain $X$. Let $Q \in \mathbb R^{N \times N}$ be the transition matrix of its associated reduced first order chain. Suppose that: 
\begin{itemize}
\item[(i)]{$Q$ has no other eigenvalues of modulus one except its dominant eigenvalue $\lambda =1$.}
\item[(ii)]{There exists $\pi \in \mathbb R^n$ such that for any (right) eigenvector $y \in \mathbb R^N$ of $Q$ corresponding to the dominant eigenvalue $\lambda$, normalized as $y \ge 0$ and $\|y\|_1=1$, we have $P^{(0)}y=\pi$ with $P^{(0)}=[\underbrace{I_n ~\ldots ~I_n}_{n^{m-2}}]$ being the mode-$1$ matricization of ${\cal P}^{(0)}={\cal I}$ and $I_n \in \mathbb R^{n \times n}$ being the identity matrix.}
\end{itemize}
Then, 
\be
\label{ro}
\lim_{k \to \infty}{\cal P}^{(k)}=\pi \otimes \underbrace{e \otimes \cdots \otimes e}_{m-1},
\ee
where $e \in \mathbb R^n$ is the vector of all ones and $\otimes$ stands for the outer product. This $\pi$ is precisely the limiting probability distribution of $X$. If, additionally, $\cal P$ is regular, then $\pi > 0$.
\end{2}

\section{Graph Theory Assisted Method}
\label{main}
\setcounter{equation}{0}

With the preparatory material in the previous sections, we are now ready to develop several new results.

It should be mentioned that if (\ref{ro}) holds, then $\pi$ is the limiting probability distribution of the higher order chain \cite{HX1}. Now, we can strengthen this by showing its converse. Before doing so, we need the following result in \cite{HX1}.

\begin{1}
Given an $(m-1)$th order, $n$-state Markov chain with transition tensor $\cal P$, let $P^{(k)}$ be the mode-$1$ matricization of ${\cal P}^{(k)}$ and let $Q$ be the transition matrix of the reduced first order chain. Then,
\be
\label{kom}
P^{(k+\ell)}=P^{(k)}Q^\ell
\ee
for any integers $k, \ell \ge 0$.
\end{1}

In the sequel, $e_i$ stands for the $i$th column of an identity matrix, whose size can be determined from the context.

\begin{2}
\label{conv}
Given an $(m-1)$th order, $n$-state Markov chain $X$ with transition tensor $\cal P$, its limiting probability distribution $\pi$ exists if and only if (\ref{ro}) holds.
\end{2}
\bp
The sufficiency part has been done in \cite{HX1}. To show necessity, we use (\ref{xt}), (\ref{yt}), and (\ref{kom}) to obtain 
$$x_{m+k}=PQ^ky_{m-1}=P^{(k+1)}y_{m-1}.$$
Hence, for any $y_{m-1}$,
\be
\label{yp}
\lim_{k \to \infty}P^{(k+1)}y_{m-1}=\pi.
\ee

Starting with $X_1=X_2=\ldots=X_{m-1}=1$, the initial probability distributions of $X$ are $x_1=x_2=\ldots=x_{m-1}=e_1$, and $y_{m-1}=e_1$. Since $\pi$ does not depend on the initial probability distributions, we know from (\ref{yp}) that the first column of $P^{(k)}$ converges to $\pi$. 

Continuing, starting with $X_1=X_2=\ldots=X_{m-2}=1$ and $X_{m-1}=2$, the initial probability distributions of $X$ are $x_1=x_2=\ldots=x_{m-2}=e_1$ and $x_{m-1}=e_2$, and $y_{m-1}=e_2$. Again from (\ref{yp}), we see that the second column of $P^{(k)}$ converges to $\pi$ too.

Proceeding in the above manner, we conclude that each column of $P^{(k)}$ converges to $\pi$, i.e., (\ref{ro}) holds.
\ep

Theorem \ref{conv} reaffirms the essential role (\ref{ro}) plays when it comes to the limiting probability distribution of a higher order Markov chain. This result also holds in the first order case. Its first order counterpart can be found in, for example, \cite{Dob}. 

To check conditions (i) and (ii) in Theorem \ref{suff}, the direct approach requires computations of the full spectrum of $Q$, plus its eigenvectors corresponding to the eigenvalue $\lambda=1$ if condition (i) is satisfied. These typically have time complexity $O(N^3)$ \cite{GV} and, therefore, are costly even for problems of moderate size. Next, we shall describe a graph theory assisted method to overcome this difficulty.

To connect a first order Markov chain with relevant results in graph theory, we start with the notion of closed communicating classes. Without loss of generality, we describe it using the reduced first order chain with transition matrix $Q=[q_{ij}]$ and state space $T=\{1, 2, \ldots, N\}$, both in the regular index form from now on. As usual, we denote $Q^k=[q^{(k)}_{ij}]$.

State $i \in T$ is said to be accessible from state $j \in T$, written as $j \rightarrow i$, if $q^{(k)}_{ij} > 0$ for some $k \ge 1$. States $i, j \in T$ are said to communicate, written as $i \leftrightarrow j$, if $i \rightarrow j$ and $j \rightarrow i$. Such a $\leftrightarrow$ relation is known to be an equivalence relation and partitions $T$ into communicating (i.e., equivalence) classes \cite{KS}. A communicating class $\alpha \subseteq T$ is called closed if $i \nrightarrow j$ whenever $i \in \alpha$ and $j \in \alpha^c$. It is also known \cite{KS} that for a first order Markov chain with a finite number of states, there exists at least one closed communicating class. 

The next lemma can be seen in, for example, \cite{Sen}, where $Q[\alpha, \beta]$ denotes the submatrix of $Q$ on rows $\alpha$ and columns $\beta$ and, by convention, $Q[\alpha]=Q[\alpha, \alpha]$. A similar notation is also used in what follows for a subvector. In addition, the spectrum and spectral radius are denoted, respectively, as $\sigma(\cdot)$ and $\rho(\cdot)$.
\begin{1}
\label{can}
With an appropriate simultaneous permutation of rows and columns if necessary, there exists integer $r \ge 1$ such that $Q \in \mathbb R^{N \times N}$ has the following canonical form: 
\be
\label{canfor}
Q=\left[\begin{array}{ccccc}
Q[\alpha_1] & 0 & \cdots & 0 & Q[\alpha_1, \beta]\\
0 & Q[\alpha_2] & \cdots & 0 & Q[\alpha_2, \beta]\\
\vdots & \vdots & & \vdots & \vdots\\
0 & 0 & \cdots & Q[\alpha_r] & Q[\alpha_r, \beta]\\
0 & 0 & \cdots & 0 & Q[\beta]
\end{array}\right],
\ee
where $\alpha_1, \alpha_2, \ldots, \alpha_r \subseteq T$ are all the closed communicating classes and $\beta=T\backslash \cup_{i=1}^r \alpha_i$. Besides, $Q[\alpha_1], Q[\alpha_2], \ldots, Q[\alpha_r]$ are irreducible column stochastic matrices, each representing a subchain consisting only of communicating states, and $Q[\beta]$ is a column substochastic matrix with $\rho(Q[\beta]) < 1$. 
\end{1}

Clearly, by Lemma \ref{can}, $\sigma(Q)=\left(\cup_{i=1}^r\sigma(Q[\alpha_i])\right) \cup \sigma(Q[\beta])$. Following the celebrated Perron-Frobenius Theorem, see \cite{Sen}, condition (i) of Theorem \ref{suff} is equivalent to the primitivity, i.e., regularity, of $Q[\alpha_1], Q[\alpha_2], \ldots, Q[\alpha_r]$. In addition, such regularity condition means that each $Q[\alpha_i]$ is irreducible and aperiodic \cite{Sen}. Besides, for each subchain that $Q[\alpha_i]$ represents, its limiting probability distribution $\pi^{(i)}$ exists. Let $y^{(i)} \in \mathbb R^N$ be such that $y^{(i)}[\alpha_i]=\pi^{(i)}$ with all the remaining entries of $y^{(i)}$ being zero. Then, it is easy to see that $y^{(i)}$ is a normalized (right) eigenvector of $Q$ corresponding to its eigenvalue $\lambda=1$. It is obvious that $y^{(1)}, y^{(2)}, \ldots, y^{(r)}$ form a basis for the eigenspace of $Q$ associated with this $\lambda$.

The above analysis yields the following counterpart of Theorem \ref{suff}.

\begin{2}
\label{irrape}
Let $Q \in \mathbb R^{N \times N}$ be the transition matrix, which is already in the canonical form (\ref{canfor}), of the reduced first order chain associated with an $(m-1)$th order, $n$-state Markov chain $X$. Suppose that:
\begin{itemize}
\item[(i)]{$Q[\alpha_1], Q[\alpha_2], \ldots, Q[\alpha_r]$ are not only irreducible but also aperiodic.}
\item[(ii)]{For each $i=1, 2, \ldots, r$, let $\pi^{(i)}$ be the normalized (right) eigenvector of $Q[\alpha_i]$ corresponding to its eigenvalue $1$, and $y^{(i)} \in \mathbb R^N$ be such that $y^{(i)}[\alpha_i]=\pi^{(i)}$ and the remaining entries being zero. Then, there exists $\pi \in \mathbb R^n$ such that $P^{(0)}y^{(i)}=\pi$ for all $i$, where $P^{(0)}=[\underbrace{I_n ~\ldots ~I_n}_{n^{m-2}}]$.}
\end{itemize}
Then, $\pi$ is the limiting probability distribution of the chain $X$. Furthermore, $\pi > 0$ if, in addition to conditions (i) and (ii) above, $X$ is regular too.
\end{2}

Here is exactly where graph theory comes in. Without loss of generality, we continue to describe the relevant concepts using $Q$ and $T$ in the regular index form. The underlying digraph $\Gamma=\Gamma(V,E)$ of $Q$ consists of the set of vertices $V=T$ and the set of directed edges $E=\{(j, i) : q_{ij}>0,  ~i, j \in V\}$, where $(j, i)$ means a directed edge from vertices $j$ to $i$. For a comprehensive reference on graph theory, see, for example, \cite{Wes}.

It is well known, see \cite{HJ}, that a nonnegative square matrix is irreducible if and only if its underlying digraph is strongly connected, i.e., there exists a directed path between every pair of vertices. Let $\Gamma[\alpha_i]$ be the underlying digraph of $Q[\alpha_i]$ in (\ref{canfor}). The irreducibility of each $Q[\alpha_i]$ translates to the strong connectedness of the respective subdigraph $\Gamma[\alpha_i]$. Moreover, $\Gamma[\alpha_i]$ is maximal in the sense that no other vertices can be added without losing its strongly connected property. Such a maximal strongly connected subdigraph is called a strongly connected component of $\Gamma$. Furthermore, since $q_{k\ell}=0$ when $\ell \in \alpha_i$ and $k \in V\backslash \alpha_i$, there are no outgoing edges leaving $\alpha_i$ and thus each $\Gamma[\alpha_i]$ is called a closed strongly connected component of $\Gamma$ \cite{Lin}. In the same vein as the closed communicating classes of a first order Markov chain with a finite number of states as mentioned earlier, there exists at least one closed strongly connected component in $\Gamma$.
 
The aperiodicity of each $Q[\alpha_i]$ can also be characterized with graph theoretic terms \cite{HJ}. A directed cycle in $\Gamma[\alpha_i]$ is a directed path with no vertex being repeated, except that it starts and ends at the same vertex. The aperiodicity of $Q[\alpha_i]$ is actually equivalent to that the greatest common divisor of the lengths of all directed cycles in the respective $\Gamma[\alpha_i]$ is $1$. In this case, the subdigraph $\Gamma[\alpha_i]$ is also called aperiodic.

Accordingly, Theorem \ref{irrape} can now be stated, in part, in graph theoretic terms as:

\begin{2}
\label{gtam}
For an $(m-1)$th order, $n$-state Markov chain $X$, let $Q \in \mathbb R^{N \times N}$ be the transition matrix of its reduced first order chain on state space $T$. Let $\Gamma[\alpha_1], \Gamma[\alpha_2], \ldots, \Gamma[\alpha_r]$ be the closed strongly connected components of $\Gamma$, the underlying digraph of $Q$, where $\alpha_1, \alpha_2, \ldots, \alpha_r \subseteq T$ (and $\cap_{i=1}^r \alpha_i=\emptyset$). Suppose that:
\begin{itemize}
\item[(i)]{$\Gamma[\alpha_1], \Gamma[\alpha_2], \ldots, \Gamma[\alpha_r]$ are all aperiodic.}
\item[(ii)]{For each $i=1, 2, \ldots, r$, let $\pi^{(i)}$ be the normalized (right) eigenvector of $Q[\alpha_i]$ corresponding to its eigenvalue $1$, and $y^{(i)} \in \mathbb R^N$ be such that $y^{(i)}[\alpha_i]=\pi^{(i)}$ and the remaining entries being zero. Then, there exists $\pi \in \mathbb R^n$ such that $P^{(0)}y^{(i)}=\pi$ for all $i$ with $P^{(0)}=[\underbrace{I_n ~\ldots ~I_n}_{n^{m-2}}]$, which is also referred to as the invariance of $P^{(0)}y^{(i)}$ with respect to $i$.}
\end{itemize}
Then, $\pi$ is the limiting probability distribution of the chain $X$. Furthermore, $\pi > 0$ if, in addition to conditions (i) and (ii) above, $X$ is also regular.
\end{2}

The above result is at the center of the graph theory assisted method for determining the limiting probability distribution for a higher order Markov chain. Specifically, this method can be described as:

\begin{7}
\label{proc}
To determine the limiting probability distribution of a higher order Markov chain, based on Theorem \ref{gtam}, we proceed as follows. 
\begin{itemize}
\item[(1)]{Construct the adjacency matrix $A$ from $Q$, the transition matrix of the reduced first order chain, and then the underlying digraph $\Gamma$ from $A$.}
\item[(2)]{Determine the strongly connected components of $\Gamma$, which can be done using the method of depth first search to index each vertex according to the order in which it is first encountered, tracking the so-called lowlink values, identifying the roots of the strongly connected components, and then forming these components \cite{Tar}.}
\item[(3)]{Examine whether each strongly connected component of $\Gamma$ is closed according to whether it has any outgoing edge. Let $\Gamma[\alpha_1], \Gamma[\alpha_2], \ldots, \Gamma[\alpha_r]$ be the resulting closed strongly connected components.}
\item[(4)]{Determine whether each $\Gamma[\alpha_i]$ is aperiodic, which is done by constructing a rooted spanning tree of $\Gamma[\alpha_i]$, with each vertex $j$ being labeled by its shortest distance $v(j)$ from the root, and then obtaining the greatest common divisor of the lengths of all directed cycles in $\Gamma[\alpha_i]$ using $v(j)-v(k)+1$ over the directed edges $(j,k)$ in $\Gamma[\alpha_i]$ \cite{Den}.}
\item[(5)]{If there is a closed strongly connected component that is not aperiodic, then condition (i) of Theorem \ref{gtam} is not satisfied and the algorithm terminates without a limiting probability distribution; otherwise, continue to the next step.}
\item[(6)]{Determine the unique normalized (right) eigenvector $\pi^{(i)}$ of each $Q[\alpha_i]$ corresponding to its eigenvalue $1$. Check condition (ii) of Theorem \ref{gtam}. If this condition is satisfied, we have $\pi$, the desired limiting probability distribution; otherwise,  the algorithm terminates without a limiting probability distribution.}
\end{itemize}
\end{7}

The main computational costs of Algorithm \ref{proc} are from steps (2), (4), and (6). The time complexity of steps (2) and (4) is $O(|V|+|E|)=O(nN)$ only \cite{Den, Tar}, where $|\cdot |$ denotes the cardinality. Step (6), on the other hand, has time complexity $O(\alpha^3)$, where $\alpha=\max_{1 \le i \le r} |\alpha_i|$. Especially, in the case that condition (i) of Theorem \ref{gtam} is not satisfied, Algorithm \ref{proc} terminates at step (5) and thus is much more effective than the direct implementation of Theorem \ref{suff}, which has time complexity $O(N^3)$ due to the computations of the full spectrum of $Q$. If condition (i) is satisfied, Algorithm \ref{proc} moves on to step (6). In this case, if $V$ is decomposed into multiple closed strongly connected components of smaller sizes, then Algorithm \ref{proc} still outperforms the direct verification of Theorem \ref{suff}. In addition, for the special case that $V$ constitutes a single closed strongly connected component, if the objective is solely to determine whether the limiting probability distribution exists or not, then again Algorithm \ref{proc} does not hinge on eigenvalue computations as required while directly checking condition (i) of Theorem \ref{suff}. Observe that for this last case, condition (ii) of Theorem \ref{gtam} is automatically satisfied. In summary, the main advantages of Algorithm \ref{proc} are twofold: (1) it employs a linear time approach \cite{Den, Tar} to validating condition (i) of Theorem \ref{suff}, and (2) it may reduce the computational costs for eigenvalues and eigenvectors, if they are needed, by carrying out such computations on submatrices of $Q$.

It is also worth mentioning here that Theorem \ref{gtam} serves as an extension to the classical result on the limiting probability distribution of a first order Markov chain as described in Section \ref{intro}, i.e., if $P$ is regular, then the limiting probability distribution exits. Observe that if the chain in Theorem \ref{gtam} is of first order, then $Q=P$. The regularity of $P$, therefore, is equivalent to that $T$ is a single closed strongly connected component, hence condition (i) there is satisfied trivially; meanwhile, condition (ii) there becomes trivial too. On the other hand, Theorem \ref{gtam} is broader because a first order chain satisfying conditions (i) and (ii) may not have a regular transition matrix; see Example \ref{ex4} in the next section.

\section{Implementation and Examples}
\label{impl}
\setcounter{equation}{0}

The graph theory assisted Algorithm \ref{proc} has been implemented using MATLAB as a function {\tt gtam} (acronym for Graph Theory Assisted Method) and added to the HOMC package as developed in \cite{Xu}. This {\tt gtam} function, along with other relevant functions such as {\tt rcmat}, {\tt eyet}, and {\tt t2mat} in the HOMC package that {\tt gtam} calls, is available at: 

\centerline{\url{https://neumann.math.siu.edu/homc}}

\noindent For descriptions of the functions {\tt rcmat}, {\tt eyet}, and {\tt t2mat}, see \cite{Xu}.

The {\tt gtam} function utilizes MATLAB's built-in functions as much as possible for better performance. Below is a brief description of the main MATLAB's built-in functions being incorporated in {\tt gtam}.

\begin{itemize}
\item[(1)]{{\tt digraph}: It constructs a digraph using its adjacency matrix, which is needed in steps (1) and (4) of Algorithm \ref{gtam}.}
\item[(2)]{{\tt conncomp}: For a digraph, it yields the strongly connected components by default, which implements step (2) of Algorithm \ref{proc}.}
\item[(3)]{{\tt distances}: It determines the shortest distances between vertices in a digraph, which is used for step (4) of Algorithm \ref{proc}.}
\item[(4)]{{\tt findedge}: It finds the directed edges in a digraph, which is also used for step (4) of Algorithm \ref{proc}.}
\item[(5)]{{\tt gcd}: It returns the greatest common divisor, which is needed in step (4) of Algorithm \ref{proc} too.}
\item[(6)]{{\tt eigs}: It calculates selected eigenvalues and eigenvectors, which helps implement step (6) of Algorithm \ref{proc}. Note that in MATLAB R2017a or earlier, its {\tt 'largestreal'} option should be replaced with {\tt 'lr'} or {\tt 'la'}.}
\end{itemize}

To use the {\tt gtam} function, its syntax goes simply as {\tt lpd=gtam(P,tol);}, where {\tt P} means the transition tensor of the higher order chain and {\tt tol} represents the tolerance while determining whether $P^{(0)}y^{(i)}$ is invariant for any $i=1, 2, \ldots, r$. The parameter {\tt tol} is optional and its default value is {\tt 1e-13}. If conditions (i) and (ii) in Theorem \ref{gtam} are both satisfied, {\tt lpd} returns with the limiting probability distribution $\pi$; otherwise, it returns with $0$.

Next, we revisit three examples from \cite{HX2} to illustrate Theorem \ref{gtam}, Algorithm \ref{proc}, and the {\tt gtam} function. In all these examples, the default value of {\tt tol} has been used.

\begin{6}
\label{ex1}
Take a second order Markov chain with transition tensor  
$${\cal P}(:,:,1)=\left[\begin{array}{ccc}
1/2 & 0 & 0\\
1/2 & 0 & 1\\
0 & 1 & 0
\end{array}\right], ~{\cal P}(:,:,2)=\left[\begin{array}{ccc}
0 & 1/2 & 1\\
0 & 1/2 & 0\\
1 & 0 & 0
\end{array}\right], ~{\cal P}(:,,3)=\left[\begin{array}{ccc}
1/2 & 1/2 & 1\\
1/2 & 1/2 & 0\\
0 & 0 & 0
\end{array}\right].$$
Using {\tt lpd=gtam(P);}, we obtain: 
\begin{verbatim}
there are   2 closed strongly connected components:
    3    4    5    8
    1    2    6    7

component   1 is aperiodic
component   2 is aperiodic

condition (i) holds

limiting probability distribution may not exist
since condition (ii) fails

\end{verbatim}
\end{6}

\begin{6}
\label{ex2}
Take a third order Markov chain with transition tensor 
$${\cal P}(:,:,1,1)=\left[\begin{array}{cc}
0 & 1/2\\
1 & 1/2
\end{array}\right], ~{\cal P}(:,:,2,1)=\left[\begin{array}{cc}
1/2 & 0\\
1/2 & 1
\end{array}\right],$$
$${\cal P}(:,:,1,2)=\left[\begin{array}{cc}
1 & 1/2\\
0 & 1/2
\end{array}\right], {\cal P}(:,:,2,2)=\left[\begin{array}{cc}
1/2 & 1\\
1/2 & 0
\end{array}\right].$$
Using {\tt lpd=gtam(P);}, we obtain: 
\begin{verbatim}
there are   1 closed strongly connected components:
    1    2    3    4    5    6    7    8

component   1 is periodic

limiting probability distribution may not exist
since condition (i) fails

\end{verbatim}
\end{6}

\begin{6}
\label{ex3}
Take a second order Markov chain with transition tensor 
$${\cal P}(:,:,1)=\left[\ba{cccc}
1/2 & 0 & 0 & 0\\
1/2 & 0 & 1 & 0\\
0 & 1 & 0 & 1\\
0 & 0 & 0 & 0
\ea\right], ~{\cal P}(:,:,2)=\left[\ba{cccc}
0 & 0 & 1/2 & 1\\
0 & 1/2 & 0 & 0\\
1/2 & 1/2 & 0 & 0\\
1/2 & 0 & 1/2 & 0
\ea\right],$$
$${\cal P}(:,:,3)=\left[\ba{cccc}
0 & 1 & 0 & 1\\
1 & 0 & 1/2 & 0\\
0 & 0 & 1/2 & 0\\
0 & 0 & 0 & 0
\ea\right], ~{\cal P}(:,:,4)=\left[\ba{cccc}
0 & 0 & 0 & 0\\
1 & 1 & 1 & 0\\
0 & 0 & 0 & 1/2\\
0 & 0 & 0 & 1/2
\ea\right].$$
Using {\tt lpd=gtam(P);}, we obtain: 
\begin{verbatim}
there are   2 closed strongly connected components:
    3    4    5   10   15
    2    7    9   12   13

component   1 is aperiodic
component   2 is aperiodic

condition (i) holds

condition (ii) holds

limiting probability distribution is:
    0.285714285714    0.285714285714    0.285714285714    0.142857142857
\end{verbatim}
We notice that the limiting probability distribution here is positive, which is due to ${\cal P}^{(10)} > 0$, i.e., the regularity of $\cal P$.
\end{6}

Finally, we give three more examples to illustrate that, as is anticipated, the graph theory assisted method applies to a first order Markov chain too. 

\begin{6}
\label{ex4}
Take a first order Markov chain with transition matrix 
$$P=\left[\begin{array}{ccc}
1 & 0 & 1\\
0 & 1/2 & 0\\
0 & 1/2 & 0
\end{array}\right].$$
Using {\tt lpd=gtam(P);}, we obtain:
\begin{verbatim}
there are   1 closed strongly connected components:
    1

component   1 is aperiodic

condition (i) holds

condition (ii) holds automatically

limiting probability distribution is:
    1.000000000000    0.000000000000    0.000000000000
\end{verbatim}
The above transition matrix $P$ is not regular, thus it is not surprising to see that the limiting probability distribution is not positive. As commented at the end of Section \ref{main}, this example also shows that conditions (i) and (ii) in Theorem \ref{gtam} are weaker, and thus more general, than the regularity requirement on $P$.
\end{6}

\begin{6}
\label{ex5}
Take a first order Markov chain with transition matrix 
$$P=\left[\begin{array}{ccc}
0 & 1/2 & 0\\
1 & 0 & 1\\
0 & 1/2 & 0
\end{array}\right].$$
Using {\tt lpd=gtam(P);}, we obtain:
\begin{verbatim}
there are   1 closed strongly connected components:
    1    2    3

component   1 is periodic

limiting probability distribution may not exist
since condition (i) fails
\end{verbatim}
\end{6}

\begin{6}
\label{ex6}
Take a first order Markov chain with transition matrix 
$$P=\left[\begin{array}{cccc}
1/3 & 0 & 1/3 & 0\\
1/3 & 1 & 0 & 0\\
1/3 & 0 & 1/3 & 0\\
0 & 0 & 1/3 & 1
\end{array}\right].$$
Using {\tt lpd=gtam(P);}, we obtain: 
\begin{verbatim}
there are   2 closed strongly connected components:
    4
    2

component   1 is aperiodic
component   2 is aperiodic

condition (i) holds

limiting probability distribution may not exist
since condition (ii) fails
\end{verbatim}
\end{6}

\section{Conclusions}
\label{conc}
\setcounter{equation}{0}

In this paper, we are primarily concerned with the verification of conditions (i) and (ii) in Theorem \ref{suff} for the existence of the limiting probability distribution of a higher order Markov chain. Directly verifying such conditions requires computations of the full spectrum of a usually very large transition matrix arising from the associated reduced first order chain and thus it may not be computationally practical. To work around this obstacle, we propose in this paper a graph theory assisted method that is formalized in Theorem \ref{gtam} and Algorithm \ref{proc}. The advantage of this method becomes particularly evident when, after verification, condition (i) is not satisfied due to the fact that its time complexity is only $O(nN)$ instead of $O(N^3)$. Even for the case that condition (i) is satisfied, if there are multiple closed strongly connected components as described in Theorem \ref{gtam},  the method still accomplishes improved time complexity. In conjection with the formalization of the method, we provide a user-friendly MATLAB implementation of the method as well, together with illustrative examples.

\end{document}